\documentstyle[twoside,12pt]{article}
\newcommand{\Proof}{{\em Proof. }}

\newtheorem{prop}{Proposition}[section]
\newtheorem{opr}[prop]{Definition}
\newtheorem{theo}[prop]{Theorem}

\newtheorem{rem}[prop]{Remark}
\newtheorem{coro}[prop]{Corollary}
\newtheorem{exam}[prop]{Example}

\date{}

\title{Quantum Cohomology Rings of
  Toric Manifolds}

\author{Victor V. Batyrev \\
 Universit\"at-GH-Essen, Fachbereich  6,  Mathematik \\
 Universit\"atsstr. 3,  45141  Essen  \\
Federal Republic of Germany \\
e-mail: matf$\emptyset \emptyset$@vm.hrz.uni-essen.de}

\begin{document}

\maketitle

\footnote{Supported by DFG, Forschungsschwerpunkt  Komplexe
Mannigfaltigkeiten and NSF (DMS-9022140).}

\begin{abstract}
We  compute the quantum cohomology
ring $H^*_{\varphi}({\bf P}_{\Sigma}, {\bf C})$
of an arbitrary  $d$-dimensio\-nal smooth projective toric
manifold ${\bf P}_{\Sigma}$ associated with a fan $\Sigma$.
The multiplicative structure of $H^*_{\varphi}({\bf P}_{\Sigma}, {\bf C})$
 depends on the choice of an element $\varphi$ in the ordinary cohomology
group $H^2({\bf P}_{\Sigma}, {\bf C})$. There are many properties
of quantum cohomology
rings $H^*_{\varphi}({\bf P}_{\Sigma}, {\bf C})$ which
are supposed to be valid for quantum cohomology rings
of  K\"ahler manifolds.
\end{abstract}

\section{Introduction}

\hspace*{\parindent}

The notion of the {\em quantum cohomology ring} of a K\"ahler manifold $V$
naturally appears in the consideration of the
so called  {\em topological sigma model}
associated with $V$ (\cite{witten}, 3a-b).
If the canonical line bundle ${\cal K}$ of
$V$ is  negative, then one  recovers the
multiplicative structure
of the quantum cohomology ring of $V$ from the intersection theory
on the moduli space ${\cal I}_{\lambda}$ of holomorphic mappings $f$ of
the complex  sphere $f\;: \; S^2 \cong {\bf CP}^1 \rightarrow V$ where
$\lambda$ is the homology class in $H_2(V, {\bf Z})$
of the image $f({\bf CP}^1)$.

If the canonical bundle ${\cal K}$ is trivial, the quantum cohomology ring
was considered by Vafa  as a important tool for explaining
the  mirror symmetry for Calabi-Yau manifolds  \cite{vafa}.

The quantum cohomology ring $QH_{\varphi}(V, {\bf C})$ of
a K\"ahler manifold $V$, unlike the ordinary cohomology ring,
have
the multiplicative structure which depends
on the class $\varphi$ of the K\"ahler
$(1,1)$-form corresponding to  a K\"ahler metric $g$ on $V$.
When we rescale the metric $g \rightarrow t g$ and
put $t \rightarrow \infty$, the quantum ring "becomes" the classical
cohomology ring. For example, for the topological sigma
model  on the
complex projective line ${\bf CP}^1$ itself,  the classical
cohomology ring is
generated by the class $x$ of a K\a"ahler $(1,1)$-form, where $x$
satisfies the quadratic equation
\begin{equation}
 x^2 = 0,
 \label{first}
\end{equation}
while the quantum  cohomology ring is also generated by $x$, but the
equation  satisfied by $x$ is different:
\begin{equation}
 x^2 = {\rm exp}(-\int_{\lambda} \varphi),
 \label{second}
\end{equation}
$\lambda$ is a  non-zero effective $2$-cycle.
Similarly, the quantum cohomology
ring of $d$-dimensional complex projective space is generated by
the element $x$ satisfying the equation
\begin{equation}
 x^{d+1} = {\rm exp}(-\int_{\lambda} \varphi ).
 \label{third}
\end{equation}

The main purpose of this paper is to construct and investigate
the quantum cohomology
ring  $QH^*_{\varphi}({\bf P}_{\Sigma},{\bf C})$
of an arbitrary smooth compact $d$-dimensional
toric manifold ${\bf P}_{\Sigma}$,
where $\varphi$ is an element of the ordinary
second cohomology group $H^2({\bf P}_{\Sigma}, {\bf C})$.
Since all projective spaces are
are  toric manifolds, we obtain a generalization of above examples of
quantum cohomology rings.

According to the physical interpretation,
a quantum cohomology ring is a closed operator algebra acting on the
fermionic Hilbert space. For example,
the equation  \ref{third}
one should better write as
an equations for the linear operator ${\cal X}$ corresponding
to the cohomology class $x$:
\begin{equation}
 {\cal X}^{d+1} = {\rm exp}(-\int_{\lambda} \varphi \,) id .
 \label{four}
\end{equation}

It is in general more convenient to define  quantum rings by polynomial
equations among generators.

\begin{opr}
{\rm Let
\[ h(t,x) = \sum_{n \in {\cal N}} c_n(t) x^n \]
be a one-parameter family of polynomials
in the polynomial ring
${\bf C}\lbrack x \rbrack$, where $x = \{x_i \}_{i \in I}$ is
a set of variables indexed by $I$, $t$ is a positive real number,
${\cal N}$ is  a fixed finite
set  of exponents.
We say that the polynomial
\[ h^{\infty}(x) = \sum_{n \in {\cal N}} c_n^{\infty} x^n \]
is the limit of $h(t,x)$ as $t \rightarrow \infty$, if
the point $\{ c_n^{\infty} \}_{n \in {\cal N}}$ of the
$(\mid {\cal N} \mid -1)$-dimensional complex projective
space is the limit
of the one-parameter family of points with homogeneous coordinates
$\{ c_n(t) \}_{n \in {\cal N}}$. }
\end{opr}

\begin{opr}
{\rm Let $R_{t}$ be a one-parameter family of commutative
algebras over {\bf C}
with a fixed set of generators $\{ r_i \}$, $t \in {\bf R}_{>0}$.
We denote by $J_{t}$ the ideal in
${\bf C} \lbrack x \rbrack$ consisting of all polynomial relations
among $\{ r_i \}$, i.e.,  the kernel of the surjective
homomorphism ${\bf C} \lbrack x \rbrack \rightarrow R_t$.
We say that the ideal $J^{\infty}$ is the limit of $J_t$ as
$t \rightarrow \infty$, if any one-parameter family of
polynomials $h(t,x) \in J_t$ has a limit, and $J^{\infty}$ is
generated as ${\bf C}$-vector space by all these limits.
The ${\bf C}$-algebra
\[ R^{\infty} = {\bf C} \lbrack x \rbrack / J^{\infty} \]
will be called the {\em limit} of $R_t$. }
\end{opr}

\begin{rem}
{\rm In general, it is not true that if
$J^{\infty} = \lim_{t \rightarrow \infty} J_t$,
and $J_t$ is generated by a finite set of polynomials
$\{ h_1(t,x) \ldots ,h_k(t,x) \}$, then $J^{\infty}$ is
generated by the limits  $\{ h_1^{\infty}(x), \ldots ,
h_k^{\infty}(x) \}$.  The limit ideal
$J^{\infty}$ is generated by the limits $h_i^{\infty}(x)$ {\em only if}
the set of polynomials $\{ h_i(t,x) \}$ form
 a Gr\"obner-type basis for $J_t$. }
\end{rem}

In this paper, we establish the following basic properties of quantum
cohomology rings of toric manifolds:
\bigskip

{\bf I} : If $\varphi$ is an element in the
interior of the K\"ahler cone $K({\bf P}_{\Sigma}) \subset
H^2({\bf P}_{\Sigma}, {\bf C})$, then there exists a  limit of
$QH^*_{t\varphi}({\bf P}_{\Sigma},{\bf C})$
as $t \rightarrow \infty$, and this limit
is isomorphic to  the ordinary cohomology ring
$H^*({\bf P}_{\Sigma},{\bf C})$ (Corollary
\ref{kon}).

{\bf II} :
Assume that two smooth projective toric manifolds
${\bf P}_{\Sigma_1}$ and ${\bf P}_{\Sigma_2}$ are isomorphic
in codimension 1, for instance, that ${\bf P}_{\Sigma_1}$
is obtained from ${\bf P}_{\Sigma_2}$ by
a flop-type birational transformation. Then the natural isomorphism
$H^2({\bf P}_{\Sigma_1}, {\bf C}) \cong
H^2({\bf P}_{\Sigma_2}, {\bf C})$ induces the isomorphism between
the quantum cohomology rings
\[  QH^*_{\varphi}({\bf P}_{\Sigma_1},{\bf C})
\cong  QH^*_{\varphi}({\bf P}_{\Sigma_2},{\bf C}) \]
(Theorem \ref{flop}). We notice that the ordinary cohomology rings of
${\bf P}_{\Sigma_1}$ and ${\bf P}_{\Sigma_2}$ are not isomorphic in general.

{\bf III} :
Assume that the first Chern class
$c_1({\bf P}_{\Sigma})$ of ${\bf P}_{\Sigma}$ belongs to the closed
K\"ahler cone $K({\bf P}_{\Sigma}) \subset H^2({\bf P}_{\Sigma}, {\bf C})$.
Then the  ring $QH^*_{\varphi}({\bf P}_{\Sigma},{\bf C})$ is
isomorphic to the Jacobian ring of a Laurent polynomial $f_{\varphi}(X)$
such that the equation $f_{\varphi}(X) = 0$ defines
an  affine Calabi-Yau hypersurface $Z_f$
in the $d$-dimensional algebraic torus $({\bf C}^*)^d$ where $Z_f$ is
"mirror symmetric" with respect to Calabi-Yau
hypersurfaces in ${\bf P}_{\Sigma}$
(Theorem \ref{mirr}).
Here by the "mirror symmetry" we mean the correspondence,
based on the polar duality {\rm \cite{bat.mir}},  between families of
Calabi-Yau hypersurfaces in toric varieties.

The properties II and III give
a general view on the recent result of
P. Aspinwall, B. Greene, and D. Morrison \cite{aspin4} who have shown,
for a family of Calabi-Yau $3$-folds $W$ that their
quantum cohomology ring $QH^*_{\varphi}(W, {\bf C})$  does not change under
a flop-type birational transformation (see also \cite{aspin1,aspin3}).
\medskip

{\bf IV}: Assume that the first Chern class
$c_1({\bf P}_{\Sigma})$ of ${\bf P}_{\Sigma}$ is divisible by
$r$, i.e., there exists a an element
$h \in H^2({\bf P}_{\Sigma}, {\bf Z})$ such that
$c_1({\bf P}_{\Sigma}) = rh$.
Then $QH^*_{\varphi}({\bf P}_{\Sigma}, {\bf C})$ has a natural
${\bf Z}_r$-grading (Theorem \ref{divis}).
We remark that $QH^*_{\varphi}({\bf P}_{\Sigma}, {\bf C})$ has no
${\bf Z}$-grading.
\medskip

The paper is organized as follows.
In Sections 2-4, we recall the definition and standard information about
toric manifolds. In Section 5, we define the quantum cohomology ring of
toric manifolds and prove their properties. In Section 6, we consider
examples of the behavior of quantum
cohomology rings under elementary birational
transformations such as blow-up and flop, we also consider the case
of singular toric varieties. In Section 7, we give an
combinatorial interpretation of the relation between the
quantum cohomology rings and the ordinary cohomology rings. In Section 8,
we show that the quantum cohomology ring can be interpreted as a Jacobian
ring of some Laurent polynomial. Finally, in Section 9, we prove
that our quantum
cohomology rings coincide with the  quantum cohomology
rings defined by $\sigma$-models on toric manifolds.
\bigskip

{\bf Acknowledgements.} It is a pleasure to acknowledge helpful discussions
with Yu. Manin, D. Morrison, Duco van Straten as well as with
S. Cecotti and C. Vafa. I would like to express my thanks for hospitality to
the Mathematical Sciences Research Institute where this work was conducted
and supported in part by the National Science Foundation (DMS-9022140), and
the DFG (Forschungsscherpunkt Komplexe Mannigfaltigkeiten).

\section{A definition of compact toric manifolds}
\hspace*{\parindent}

The {\em  toric varieties} were  considered in full generality in
\cite{dan2,oda1}. For the general definition
of toric variety which includes affine
and quasi-projective toric varieties with  singularities,  it is more
convenient to use the language of schemes. However,  for our purposes,
it will be sufficient to have a simplified more classical
version of the definition  for
{\em smooth and compact toric varieties}
over the complex number field ${\bf C}$. In fact, this approach to
compact toric manifolds was first proposed by M. Audin \cite{audin},
and developed by D. Cox \cite{cox}.
\bigskip

In order to obtain a $d$-dimensional compact toric manifold $V$, we need
a combinatorial object $\Sigma$, a {\em complete fan of regular
cones},  in a $d$-dimensional vector space over
${\bf R}$.
\bigskip

Let $N$, $M = {\rm Hom}\, (N, {\bf Z})$ be dual lattices of rank $d$, and
$N_{\bf R}$, $M_{\bf R}$ their ${\bf R}$-scalar extensions to $d$-dimensional
vector spaces.

\begin{opr}
{\rm  A convex subset $\sigma \subset N_{\bf R}$ is called a {\em regular
$k$-dimensional  cone} $(k \geq 1)$
if there exist $k$ linearly independent
elements $v_1, \ldots, v_k \in N$ such that
\[ \sigma = \{ \mu_1 v_1 + \cdots + \mu_k v_k \mid
\mu_i \in {\bf R}, \mu_i \geq 0 \},  \]
and the set $\{ v_1, \ldots, v_k \}$ is a subset of some ${\bf Z}$-basis of
$N$.
In this case, we call $v_1, \ldots, v_k \in N$ the
{\em integral generators of} $\sigma$.

The origin $0 \in N_{\bf R}$ we call the {\em regular
$0$-dimensional  cone}. By definition, the set of
integral generators of this cone is empty.}
\end{opr}

\begin{opr}
{\rm A regular  cone $\sigma'$ is called {\em a face} of a
regular  cone $\sigma$ (we write $\sigma' \prec \sigma)$ if the
set of integral generators of $\sigma'$ is a subset of the
set of  integral generators of $\sigma$. }
\label{face}
\end{opr}

\begin{opr}
{\rm A finite system $\Sigma = \{ \sigma_1, \ldots , \sigma_s \}$ of
regular cones in $N_{\bf R}$ is called
{\em a complete  d-dimensional
fan} of regular cones,  if the following conditions are satisfied:

(i) if $\sigma \in  \Sigma$ and $\sigma' \prec \sigma$, then
$\sigma' \in \Sigma$;

(ii) if $\sigma$, $\sigma'$ are in $\Sigma$, then
$\sigma \cap \sigma' \prec \sigma$ and $\sigma \cap \sigma' \prec \sigma'$;

(iii) $N_{\bf R} = \sigma_1 \cup \cdots \cup \sigma_s$. \\
The set of all $k$-dimensional cones in $\Sigma$ will be denoted by
$\Sigma^{(k)}$. }
\label{def.fan}
\end{opr}

\begin{exam}
{\rm  Choose $d+1$ vectors
$ v_1, \ldots, v_{d+1}$ in a
$d$-dimensional real space $E$ such that $E$ is spanned by
$v_1, \ldots, v_{d+1}$
and there exists the linear relation
\[  v_1 + \cdots +  v_{d+1} =0. \]
Define $N$ to be the lattice in $E$ consisting of all integral linear
combinations of $v_1, \ldots, v_{d+1}$. Obviously,
$N_{\bf R} = E$. Then any $k$-element subset $I \subset
\{ v_1, \ldots, v_{d+1} \}$ $(k \leq d)$ generates a
$k$-dimensional regular cone $\sigma(I)$.
The set $\Sigma(d)$ consisting of $2^{d+1} -1$ cones
$\sigma(I)$ generated by $I$ is
 a complete $d$-dimensional fan of regular cones. }
\label{weig}
\end{exam}

\begin{opr}
{\rm (cf. \cite{bat.class})
Let $\Sigma$ be a complete $d$-dimensional fan of regular cones. Denote by
$G(\Sigma) =
\{ v_1, \ldots , v_n \}$ the set of all generators of $1$-dimensional cones
in $\Sigma$ ($n = {\rm Card}\, \Sigma^{(1)}$). We call a subset
${\cal P} =\{ v_{i_1}, \ldots , v_{i_p} \} \subset G(\Sigma)$ a
{\em primitive collection} if $\{v_{i_1}, \ldots , v_{i_p}\}$ is not
the set of generators of a $p$-dimensional simplicial cone in $\Sigma$,
while for all $k$ $(0 \leq k < p)$ each $k$-element subset of
${\cal P}$ generates a $k$-dimensional cone in $\Sigma$.}
\end{opr}

\begin{exam}
{\rm Let $\Sigma$ be a fan $\Sigma(d)$ from Example \ref{weig}. Then
there exists the unique  primitive collection ${\cal P} =
G(\Sigma(d))$.}
\label{ex.prim}
\end{exam}

\begin{opr}
{\rm Let ${\bf C}^n$
be  $n$-dimensional affine space over ${\bf C}$ with the set of coordinates
$z_1, \ldots, z_n$ which are in the one-to-one correspondence
$z_i \leftrightarrow v_i$ with elements of $G(\Sigma)$.
Let ${\cal P} =\{ v_{i_1}, \ldots , v_{i_p} \}$ be a primitive collection
in $G(\Sigma)$. Denote by  ${\bf A}({\cal P})$  the $(n-p)$-dimensional
affine subspace in ${\bf C}^n$  defined by  the equations
\[ z_{i_1} = \cdots = z_{i_p} = 0. \]}
\end{opr}

\begin{rem}
{\rm Since every primitive collection ${\cal P}$ has at least two elements,
the codimension of ${\bf A}({\cal P})$ is at least $2$.}
\label{cod.2}
\end{rem}

\begin{opr}
{\rm Define the closed algebraic subset $Z(\Sigma)$
in ${\bf C}^n$ as follows
\[ Z(\Sigma) = \bigcup_{\cal P} {\bf A}({\cal P}), \]
where ${\cal P}$ runs over all primitive collections in $G(\Sigma)$.
Put }
\[ U(\Sigma) = {\bf C}^n \setminus Z(\Sigma). \]
\end{opr}

\begin{opr}
{\rm Two complete $d$-dimensional fans of regular cones  $\Sigma$ and
$\Sigma'$ are called {\em combinatorially equivalent} if there exists a
bijective mapping $\Sigma \rightarrow \Sigma'$ respecting the
face-relation "$\prec$" (see \ref{face}).}
\end{opr}

\begin{rem}
{\rm It is easy to see that the open subset
$U(\Sigma) \subset {\bf C}^n$ depends only on the
combinatorial structure of $\Sigma$, i.e.,  for any
two combinatorially equivalent
fans $\Sigma$ and $\Sigma'$, one has $U(\Sigma) \cong U(\Sigma')$. }
\end{rem}

\begin{opr}
{\rm Let $R(\Sigma)$ be the subgroup  in ${\bf Z}^n$ consisting of all
lattice vectors $\lambda = (\lambda_1, \ldots, \lambda_n)$ such that
$\lambda_1 v_1 + \cdots + \lambda_n v_n = 0$. }
\end{opr}

Obvioulsy, $R(\Sigma)$ is isomorphic to ${\bf Z}^{n-d}$.

\begin{opr}
{\rm Let $\Sigma$ be a complete $d$-dimensional fan of regular cones.
Define ${\bf D}(\Sigma)$ to be the connected commutative
subgroup in $({\bf C}^*)^n$ generated by all one-parameter subgroups
\[ a_{\lambda}\;\; : \;\; {\bf C}^* \rightarrow ({\bf C}^*)^n, \;\; \]
\[ t \rightarrow (t^{\lambda_1}, \ldots, t^{\lambda_n}) \]
where $\lambda \in R(\Sigma)$.}
\end{opr}

\begin{rem}
{\rm Choosing a {\bf Z}-basis in $R(\Sigma)$, one easily obtains an
isomorphism between ${\bf D}(\Sigma)$ and $({\bf C}^*)^{n-d}$. }
\end{rem}
\bigskip

Now we are ready to give the definition of the compact toric manifold
${\bf P}_{\Sigma}$ associated with a complete $d$-dimensional
fan of regular cones $\Sigma$.
\bigskip

\begin{opr}
{\rm Let $\Sigma$ be a complete $d$-dimensional fan of regular
cones.
Then quotient
\[ {\bf P}_{\Sigma} = U(\Sigma)/{\bf D}(\Sigma) \]
is called the {\em  compact
toric manifold associated with} $\Sigma$. }
\end{opr}

\begin{exam}
{\rm  Let $\Sigma$ be a fan $\Sigma(d)$ from Example \ref{weig}.
By \ref{ex.prim},
$U(\Sigma(d)) = {\bf C}^{d+1} \setminus \{0 \}$. By the definition
 of $\Sigma(d)$, the subgroup $R(\Sigma) \subset {\bf Z}^n$ is generated
 by $(1, \ldots, 1) \in {\bf Z}^{d+1}$.
 Thus, ${\bf D}(\Sigma) \subset ({\bf C}^*)^n$  consists of
the elements  $(t, \ldots , t)$, where $t \in {\bf C}^*$.
So the toric manifold
 associated with $\Sigma(d)$ is the ordinary $d$-dimensional projective
 space.}
\end{exam}

A priori, it is not obvious that the quotient space
${\bf P}_{\Sigma} = U(\Sigma)/{\bf D}(\Sigma)$ always exists as the
space of orbits  of the group ${\bf D}(\Sigma)$ acting free on $U(\Sigma)$,
and that ${\bf P}_{\Sigma}$ is smooth and compact. However, these facts
are easy to check if we take the $d$-dimensional projective space
${\bf P}_{\Sigma(d)}$ as a model example.
\bigskip

There exists a simple open covering of $U(\Sigma)$ by affine algebraic
varieties:

\begin{prop}
Let $\sigma$ be a $k$-dimensional cone in $\Sigma$ generated by
$\{ v_{i_1}, \ldots, v_{i_k} \}$. Define the open subset
$U(\sigma) \subset {\bf C}^n$ by the conditions
$ z_{j} \neq 0\;\; {\rm for\;\; all}\;\; j \not\in \{ i_1, \ldots ,
 i_k \}$.
 Then the open subsets $U(\sigma)$ $(\sigma \in \Sigma)$
have the properties:

{\rm (i)}
\[ U(\Sigma) = \bigcup_{\sigma \in \Sigma} U(\sigma);\]

{\rm (ii)} if $\sigma \prec \sigma'$, then $U(\sigma) \subset U(\sigma')$$;$

{\rm (iii)} for any two cone $\sigma_1, \sigma_2 \in \Sigma$, one
has
$ U(\sigma_1) \cap U(\sigma_2) = U(\sigma_1 \cap \sigma_2)$;
in particular,
\[U(\Sigma) = \bigcup_{\sigma \in \Sigma^{(d)}} U(\sigma). \]
 \end{prop}

\begin{prop}
Let $\sigma$ be a $d$-dimensional cone in $\Sigma^{(d)}$ generated by
$\{ v_{i_1}, \ldots, v_{i_d} \} \subset N$. Denote by
$ u_{i_1}, \ldots, u_{i_d}$ the dual to $ v_{i_1}, \ldots, v_{i_d}$
${\bf Z}$-basis of the  lattice $M$, i.e,
$\langle v_{i_k}, u_{i_l} \rangle = \delta_{k,l}$,
where
$\langle * , *  \rangle\;\;:\;\; N \times M \rightarrow {\bf Z}$
is the canonical pairing between lattices $N$ and $M$.

Then the affine open subset $U(\sigma)$ is isomorphic to ${\bf C}^d \times
({\bf C}^*)^{n-d}$, the action of ${\bf D}(\Sigma)$ on $U(\sigma)$
is free, and the space of ${\bf D}(\Sigma)$-orbits is isomorphic to the
affine space $U_{\sigma} = {\bf C}^d$ whose coordinate
functions $x_1^{\sigma}, \ldots, x_d^{\sigma}$ are $d$ Laurent
monomials in $z_1, \ldots, z_n$ $:$
\[ x_1^{\sigma} = z_1^{\langle v_1,u_{i_1} \rangle} \cdots ,
z_n^{\langle v_n,u_{i_1} \rangle}, \; \ldots \;,
x_d^{\sigma} = z_1^{\langle v_1,u_{i_d} \rangle}  \cdots
z_n^{\langle v_n,u_{i_d} \rangle}. \]
\end{prop}

The last statement yields a general formula for  the local
affine coordinates $x_1^{\sigma}, \ldots, x_d^{\sigma}$ of a point
$p \in U_{\sigma}$
as functions of its "homogeneous coordinates" $z_1, \ldots, z_n$ (see also
\cite{cox}).

Compactness of ${\bf P}_{\Sigma}$ follows from the fact that the
local polydiscs
\[ D_{\sigma} = \{ x \in U_{\sigma}\; : \;
\mid x_1^{\sigma} \mid \leq 1, \dots, \mid x_d^{\sigma} \mid \leq 1 \},
\; \; \sigma \in \Sigma^{(d)}  \]
form a finite compact covering of ${\bf P}_{\Sigma}$.

\section{Cohomology of toric manifolds}

\hspace*{\parindent}

Let $\Sigma$ be a complete $d$-dimensional fan of regular cones.

\begin{opr}
{\rm A continious function
$\varphi\; :\; N_{\bf R} \rightarrow {\bf R}$
 is called $\Sigma$-{\em piecewise linear}, if $\sigma$ is a linear
 function on every cone $\sigma \in \Sigma$.}
\end{opr}

\begin{rem}
{\rm It is clear that any $\Sigma$-piecewise linear function $\varphi$
is uniquely defined by its values on elements $v_i$ of $G(\Sigma)$.
So the space of
all $\Sigma$-piecewise linear functions $PL(\Sigma)$ is canonically
isomorphic to ${\bf R}^n$:
$ \varphi \rightarrow (\varphi(v_1), \ldots, \varphi(v_n)).$}
\end{rem}

\begin{theo}
The  space $PL(\Sigma)/M_{\bf R}$
of all $\Sigma$\--piece\-wise li\-near func\-tions mo\-dulo
the $d$-dimen\-sional sub\-space of glo\-bally li\-near func\-tions on
$N_{\bf R}$ is
canonically isomorphic to the cohomology space
$ H^{2}({\bf P}_{\Sigma}, {\bf R})$.
Moreover,
the first Chern class $c_1({\bf P}_{\Sigma})$, as an element of
$H^{2}({\bf P}_{\Sigma}, {\bf Z})$, is represented by the class of
the ${\Sigma}$-piecewise linear function
$\alpha_{\Sigma} \in PL(\Sigma)$ such that $\alpha_{\Sigma}(v_1) = \cdots =
\alpha_{\Sigma}(v_n) = 1$.
\end{theo}

\begin{theo}
Let $R(\Sigma)_{\bf R}$ be the ${\bf R}$-scalar extension of the abelian
group $R(\Sigma)$. Then the space $R(\Sigma)_{\bf R}$ is
canonically isomorphic to the homology space
$H_{2}({\bf P}_{\Sigma}, {\bf R})$.
\end{theo}

\begin{opr}
{\rm Let $\varphi$ be an element of $PL(\Sigma)$, $\lambda$ an element of
$R(\Sigma)_{\bf R}$.
Define the {\em degree of $\lambda$ relative to $\varphi$} as
\[ {\rm deg}_{\varphi} (\lambda) = \sum_{i =1}^n \lambda_i \varphi(v_i). \]
 }
\label{deg}
\end{opr}

It is easy to see that
for any $\varphi \in M_{\bf R}$ and for any  $\lambda \in
R(\Sigma)_{\bf R}$,
one has ${\rm deg}_{\varphi} (\lambda) = 0$. Moreover,
the degree-mapping induces a nondegenerate
pairing
\[ {\rm deg}\; :\;
PL(\Sigma)/M_{\bf R} \times R(\Sigma)_{\bf R} \rightarrow {\bf R} \]
which coincides with the canonical intersection pairing
\[H^{2}({\bf P}_{\Sigma}, {\bf R}) \times
H_2({\bf P}_{\Sigma}, {\bf R}) \rightarrow {\bf R}. \]

\begin{opr}
{\rm
Let ${\bf C}\lbrack z \rbrack$ be the polynomial ring in $n$ variables
$z_1, \ldots, z_n$.
Denote by  $SR(\Sigma)$ the ideal in
${\bf C} \lbrack z \rbrack$
generated by all monomials
\[ \prod_{v_j \in {\cal P}} z_j, \]
where ${\cal P}$ runs over all primitive collections in $G(\Sigma)$.
The ideal  $SR(\Sigma)$ is usually called the {\em Stenley-Reisner
ideal} of ${\Sigma}$.
}
\end{opr}

\begin{opr}
{\rm Let $u_1, \ldots, u_d$ be any ${\bf Z}$-basis of the lattice $M$.
Denote by
$P(\Sigma)$ the ideal in ${\bf C} \lbrack z \rbrack$ generated by $d$
elements
\[ \sum_{i = 1}^{n} \langle v_i , u_1 \rangle z_i, \ldots ,
 \sum_{i = 1}^{n} \langle v_i , u_d \rangle z_i. \]}
\end{opr}

Obviously, the ideal $P(\Sigma)$ does not
depend on the choice of basis of
$M$.

\begin{theo}
The cohomology ring of the compact toric manifold ${\bf P}_{\Sigma}$ is
canonically isomorphic to the quotient of
${\bf C} \lbrack z \rbrack$ by the sum of two ideals $P(\Sigma)$ and
$SR(\Sigma)$:
\[ H^*({\bf P}_{\Sigma}, {\bf C}) \cong {\bf C} \lbrack z \rbrack /
(P(\Sigma) + SR(\Sigma)). \]
Moreover, the canonical embedding
$ H^2({\bf P}_{\Sigma}, {\bf C}) \hookrightarrow
H^*({\bf P}_{\Sigma}, {\bf C})$
is induced by the linear mapping
\[ PL(\Sigma)\otimes_{\bf R} {\bf C}
\rightarrow {\bf C}\lbrack z \rbrack, \;\;
 \varphi \mapsto \sum_{i =1}^n \varphi_i(v_i) z_i. \]
In particular, the first Chern class of ${\bf P}_{\Sigma}$ is represented
by the sum $ z_1 + \cdots + z_n$.
\label{ord}
\end{theo}

\begin{exam}
{\rm Let ${\bf P}_{\Sigma}$ be $d$-dimensional projective space defined
by the fan $\Sigma(d)$ (see \ref{weig}). Then
\[ P(\Sigma(d)) = < (z_1 - z_{d+1}), \ldots, (z_d - z_{d+1}) >, \]
\[ SR(\Sigma(d)) = < \prod_{i=1}^{d+1} z_i  >. \]
So we obtain
\[ {\bf C} \lbrack z_1, \ldots, z_{d+1} \rbrack /
(P(\Sigma(d)) + SR(\Sigma(d))  \cong {\bf C} \lbrack x \rbrack / x^{d+1}. \]
}
\end{exam}
\bigskip

\section{Line bundles and K\"ahler classes}
\hspace*{\parindent}

Let
\[ \pi \; : \; U({\Sigma}) \rightarrow {\bf P}_{\Sigma} \]
be the canonical projection whose fibers are principal
homogeneous spaces of ${\bf D}(\Sigma)$. For any line bundle ${\cal L}$
over ${\bf P}_{\Sigma}$, the pullback $\pi^* {\cal L}$ is a line
bundle over $U(\Sigma)$. By \ref{cod.2}, $\pi^* {\cal L}$ is isomorphic
to ${\cal O}_{U(\Sigma)}$.
Therefore, the Picard group of ${\bf P}_{\Sigma}$ is isomorphic to
the group of all ${\bf D}$-linearization of ${\cal O}_{U(\Sigma)}$,
or to the group of all characters
\[ \chi\; : \; {\bf D}(\Sigma) \rightarrow {\bf C}^*.   \]
The latter is isomorphic to the group ${\bf Z}^n /M$ where ${\bf Z}^n$ is
the group of all ${\Sigma}$-piecewise linear functions $\varphi$ such that
$\varphi (N) \subset {\bf Z}$.

\begin{prop}
Assume that a character $\chi$ is represented by the class of
an integral  ${\Sigma}$-piecewise linear function $\varphi$. Then
the space
\[ H^0 ({\bf P}_{\Sigma}, {\cal L}_{\chi} ) \]
of global sections  of the corresponding line bundle ${\cal L}_{\chi}$,
is canonically isomorphic to the space of all polynomials
$F(z_1, \ldots, z_n ) \in {\bf C} \lbrack z \rbrack$ satisfying the
condition
\[ F(t^{\lambda_1}z_1, \ldots, t^{\lambda_n}z_n ) =
t^{{\rm deg}_{\varphi} \lambda} F(z_1, \ldots, z_n ) \]
\[ {for \;\; all}\;\; \lambda \in R(\Sigma),\; t \in {\bf C}^*. \]

The exponentials $(m_1, \ldots, m_n)$ of the monomials satisfying
the above condition can be identified with integral points
in the convex  polyhedron:
\[ \Delta_{\varphi} = \{ (x_1, \ldots, x_n) \in {\bf R}^n_{\geq 0} :
{\rm deg}_{\varphi} \lambda = \lambda_1 x_1 + \cdots + \lambda_n x_n,
\;\;  \lambda \in R(\Sigma) \} \]
\label{poly1}
 \end{prop}

\begin{opr}
{\rm A $\Sigma$-piecewise linear function
$\varphi \in PL(\Sigma)$
is called a {\em strictly convex support} function for the
fan $\Sigma$, if $\varphi$ satisfies the properties

(i) $\varphi$ is an {\em upper convex function}, i.e.,
\[ \varphi(x) + \varphi(y) \geq \varphi(x +y ); \]

(ii) for any two different $d$-dimensional cones $\sigma_1$, $\sigma_2
\in \Sigma$, the restrictions
$\varphi {\mid}_{\sigma}$  and
$ \varphi {\mid}_{\sigma ' }$   are different linear functions. }
\label{support}
\end{opr}

\begin{prop}
If $\varphi$ is a strictly convex support function, then the polyhedron
$\Delta_{\varphi}$ is simple $($ i.e., any vertex of $\Delta_{\varphi}$ is
contained in $d$-faces of codimension $1$$)$, and the fan $\Sigma$ can be
uniquely recovered from $\Delta_{\varphi}$ using the property:
\[ \Delta_{\varphi} \cong \{ x \in M_{\bf R}\; : \;
\langle v_i, x \rangle \geq -\varphi({v_i}) \}. \]
\end{prop}

\begin{opr}
{\rm Denote by  $K(\Sigma)$ the cone in  $H^{2}({\bf P}_{\Sigma}, {\bf R}) =
PL(\Sigma)/M_{\bf R}$ consisting
of the classes  of all upper convex ${\Sigma}$-piecewise linear
functions $\varphi \in PL(\Sigma)$.
We denote by $K^0(\Sigma)$ the interior of $K(\Sigma)$,
i.e., the cone consisting of the classes of all  strictly convex
support functions in $PL(\Sigma)$.}
\end{opr}

\begin{theo}
The open cone $K^0(\Sigma) \subset H^{2}({\bf P}_{\Sigma}, {\bf R})$
consists of classes of K\"ahler $(1,1)$-forms on ${\bf P}_{\Sigma}$,
i.e., $K(\Sigma)$ is isomorphic to the closed K\"ahler cone of
${\bf P}_{\Sigma}$.
\end{theo}
\bigskip

Next theorem will play the central role in the sequel. Its statement
is contained implicitly in \cite{oda.park,reid}:
\bigskip

\begin{theo}
A $\Sigma$-piecewise linear function $\varphi$ is a strictly convex
support function, i.e., $\varphi \in K^0(\Sigma)$,   if and only if
\[ \varphi(v_{i_1}) + \cdots + \varphi(v_{i_k}) >
\varphi(v_{i_1} + \cdots + v_{i_k}) \]
for all primitive collections ${\cal P} = \{ v_{i_1}, \ldots ,
v_{i_k} \}$ in $G(\Sigma)$.
\label{crit}
\end{theo}

\section{Quantum cohomology rings}
\hspace*{\parindent}

\begin{opr}
{\rm Let $\varphi$ be a ${\Sigma}$-piecewise linear function with complex
values, or an element  of the complexified space $PL(\Sigma)_{\bf C} =
PL(\Sigma)\otimes_{\bf R} {\bf C}$. Define the quantum cohomology ring as
the quotient of the polynomial ring
${\bf C} \lbrack z  \rbrack$ by the sum of ideals
$P(\Sigma)$ and $Q_{\varphi}(\Sigma)$:
\[ QH^*_{\varphi}({\bf P}_{\Sigma}, {\bf C}) : =
{\bf C} \lbrack z  \rbrack/ (P(\Sigma) + Q_{\varphi}(\Sigma)) \]
where
$Q_{\varphi}(\Sigma)$ is generated by binomials
\[ \exp ( \sum_{i = 1}^n a_i \varphi(v_i))
\prod_{i =1}^n z_j^{a_j} -
\exp ( \sum_{j = 1}^n b_j \varphi(v_j)) \prod_{j =1}^n z_j^{b_j} \]
running over all possible linear relations
\[ \sum_{i = 1}^n a_i v_i =  \sum_{j = 1}^n b_j v_j , \]
where all coefficients $a_i$ and $b_j$ are non-negative and integral. }
\label{def.quant}
\end{opr}

\begin{opr}
{\rm Let ${\cal P}  = \{ v_{i_1}, \ldots, v_{i_k} \} \subset G(\Sigma)$
be a primitive collection, $\sigma_{\cal P}$ the
minimal cone in $\Sigma$ containing
the sum
\[  v_{\cal P} = v_{i_1} + \ldots + v_{i_k}, \]
$v_{j_1}, \ldots,  v_{j_l}$ generators of $\sigma_{\cal P}$. Let
$l$ be the dimension of $\sigma_{\cal P}$.
By \ref{def.fan}(iii), there exists the unique representation
of $v_{\cal P}$
as an integral
linear combination of generators
$v_{j_1}, \ldots,  v_{j_l}$ with positive integral coefficients
$c_1, \ldots, c_l$:
\[ v_{\cal P} =  c_1 v_{j_1} + \cdots + c_l v_{j_l}, \]
We put
\[ E_{\varphi}({\cal P}) = \exp ( \varphi(v_{i_1} + \ldots + v_{i_k})
- \varphi(v_{i_1}) - \ldots - \varphi(v_{i_k})) \]
\[ = \exp (  c_1 \varphi(v_{j_1}) + \cdots + c_l \varphi(v_{j_l})
- \varphi(v_{i_1}) - \ldots - \varphi(v_{i_k})).\]}
\end{opr}

\begin{theo}
Assume that the K\"ahler cone $K(\Sigma)$ has the non-empty interior, i.e.,
${\bf P}_{\Sigma}$ is projective. Then the
ideal $Q_{\varphi}(\Sigma)$ is generated by
the binomials
\[ B_{\varphi}({\cal P}) = z_{i_1}  \cdots  z_{i_k} - E_{\varphi}({\cal P})
z_{j_1}^{c_1}  \cdots z_{j_l}^{c_l}, \]
where ${\cal P}$ runs over all primitive collections in $G(\Sigma)$.
\label{basis}
\end{theo}

\Proof
We use some ideas from \cite{strum}.
Let $\phi$ be an element in $PL(\Sigma)$ representing
an interior point of $K(\Sigma)$. Define the weights $\omega_1,
\ldots, \omega_n$ of $z_1, \ldots, z_n$ as
\[ \omega_i = \phi(v_i)\;\; (1 \leq i \leq n).\]
We claim that binomials $B_{\varphi}({\cal P})$ form a reduced
Gr\"obner basis for $Q_{\varphi}(\Sigma)$ relative to the
weight vector
\[ \omega = (\omega_1, \ldots, \omega_n ).\]
Notice that the weight of the monomial $z_{i_1}  \cdots  z_{i_k}$
is greater than the weight of the monomial
$z_{j_1}^{c_1}  \cdots z_{j_l}^{c_l}$, because
\[ \phi(v_{i_1}) +  \cdots + \phi(v_{i_k}) >
\phi(v_{i_1} +  \cdots + v_{i_k}) =c_1  \phi(v_{j_1}) + \cdots
c_l \phi(v_{j_l}) \]
(Theorem \ref{crit}). So the initial ideal
$init_{\omega} \langle B_{\varphi}({\cal P}) \rangle$ of the ideal
$\langle B_{\varphi}({\cal P}) \rangle$ generated by
$B_{\varphi}({\cal P})$ coincides with the ideal $SR(\Sigma)$.
It suffices to
show that the initial ideal ${\em init}_{\omega} Q_{\varphi}(\Sigma)$
also equals
$SR(\Sigma)$. The latter again follows
from Theorem \ref{crit}. \hfill $\Box$
\bigskip

\begin{opr}
{\rm The tube domain  in the complex cohomology space
$H^2({\bf P}_{\Sigma}, {\bf C})$:
\[ K(\Sigma)_{\bf C} = K(\Sigma) + i H^2({\bf P}_{\Sigma}, {\bf R}) \]
we call the {\em complexified K\"ahler cone of } ${\bf P}_{\Sigma}$. }
\end{opr}

\begin{coro}
Let $\varphi$ be an element of $H^2({\bf P}_{\Sigma}, {\bf C})$, $t$ a
positive real number. Then all
generators
$B_{t\varphi}({\cal P})$ of the ideal $Q_{t\varphi}(\Sigma)$ have finite
limits as $t \rightarrow \infty$ if and only if $\varphi \in
K(\Sigma)_{\bf C}$. Moreover, if $\varphi \in
K(\Sigma)_{\bf C}$, then the limit of
$ QH^*_{t\varphi}({\bf P}_{\Sigma}, {\bf C}) $
is  the ordinary cohomology ring
$H^*({\bf P}_{\Sigma}, {\bf C})$.
\label{kon}
\end{coro}

\Proof  Applying  Theorem \ref{crit}, we obtain:
\[ \lim_{t \rightarrow \infty} B_{t\varphi}({\cal P}) =
z_{i_1} \cdots z_{i_k}. \]
Thus,
\[ \lim_{t \rightarrow \infty} Q_{t\varphi}(\Sigma) = SR(\Sigma). \]
By Theorem \ref{ord},
\[ \lim_{t \rightarrow \infty}  QH^*_{t\varphi}({\bf P}_{\Sigma}, {\bf C})
= H^*({\bf P}_{\Sigma}, {\bf C}). \]
\hfill $\Box$
\bigskip

\begin{exam}
{\rm Consider the fan $\Sigma(d)$ defining $d$-dimensional projective space
(see \ref{weig}). Then we obtain
\[ QH^*_{\varphi}({\bf P}_{\Sigma},
{\bf C}) \cong {\bf C} \lbrack x \rbrack /
(x^{d+1} - \exp ( - {\rm deg}_{ \varphi} \lambda  )), \]
where $\lambda = (1,\ldots, 1)$ is the generator of $R(\Sigma(d))$.
This  shows the quantum cohomology ring
$QH^*_{\varphi}({\bf CP}^d, {\bf C})$
coincides with the quantum cohomology ring for ${\bf CP}^d$ proposed by
physicists. }
\end{exam}

It is important to remark that the quantum cohomology ring
$QH^*_{\varphi}({\bf P}_{\Sigma}, {\bf C})$ has no any
${\bf Z}$-grading, but
it is possible to define a ${\bf Z}_N$-grading on it.

\begin{theo}
{\rm Assume that the first Chern class $c_1({\bf P}_{\Sigma})$
is divisible by $r$. Then the ring
$QH^*_{\varphi}({\bf P}_{\Sigma}, {\bf C})$ has a
natural ${\bf Z}_r$-grading.}
\label{divis}
\end{theo}

\Proof  A linear relation
\[ \sum_{i = 1}^n a_i v_i = \sum_{j = 1}^n b_j v_j \]
gives rise to an element
\[ \lambda = (a_1 - b_1, \ldots, a_n - b_n ) \in R(\Sigma). \]
By our assumption,
\[ {\rm deg}_{\alpha_{\Sigma}} \lambda =
\sum_{i = 1}^n a_i - \sum_{j = 1}^n b_j \]
is the intersection number of $c_1({\bf P}_{\Sigma})$ and $\lambda \in
H_2 ({\bf P}_{\Sigma}, {\bf C})$, i.e., it
is divisible by $r$. This means that the binomials
\[ \exp (\sum_{i = 1}^n a_i \varphi(v_i))
\prod_{i =1}^n z_j^{a_j} -
\exp (\sum_{j = 1}^n b_j \varphi(v_j)) \prod_{j =1}^n z_j^{b_j} \]
are ${\bf Z}_r$-homogeneous. \hfill $\Box$
\bigskip

Although, the quantum cohomology ring $QH^*_{\varphi}({\bf P}_{\Sigma},
{\bf C})$ has no any ${\bf Z}$-grading, it is possible to define a
graded version of
this quantum cohomology ring over the Laurent polynomial ring
${\bf C} \lbrack z_0, z_0^{-1} \rbrack$.

\begin{opr}
{\rm Let $\varphi$ be a ${\Sigma}$-piecewise linear function with complex
values from the complexified space $PL(\Sigma)_{\bf C} =
PL(\Sigma)\otimes_{\bf R} {\bf C}$. Define the quantum cohomology ring
\[ QH^*_{\varphi}
({\bf P}_{\Sigma}, {\bf C} \lbrack z_0, z_0^{-1} \rbrack) \]
as
the quotient of the Laurent polynomial extension
${\bf C} \lbrack z  \rbrack \lbrack z_0, z_0^{-1} \rbrack$
by the sum of ideals
$Q_{\varphi,z_0}(\Sigma)$ and $P(\Sigma)$:
where
$Q_{\varphi,z_0}(\Sigma)$ is generated by binomials
\[ \exp ( \sum_{i = 1}^n a_i \varphi(v_i))
z_0^{(-\sum_{i = 1}^n a_i)}
 \prod_{i =1}^n z_i^{a_i} -
\exp ( \sum_{j = 1}^n b_j \varphi(v_j)) z_0^{(-
\sum_{j = 1}^n b_j)}
\prod_{l =1}^n z_j^{b_j} \]
running over all possible linear relations
\[ \sum_{i = 1}^n a_i v_i = \sum_{j = 1}^n b_j v_j \]
with non-negative integer coefficients $a_i$ and $b_j$. }
\label{def.quant1}
\end{opr}

The properties of the ${\bf Z}$-graded quantum cohomology ring
\[ QH^*_{\varphi}({\bf P}_{\Sigma}, {\bf C}
\lbrack z_0, z_0^{-1} \rbrack) \]
 are analogous to the properties of
$QH^*({\bf P}_{\Sigma}, {\bf C})$:

\begin{theo}
For every binomial  $B_{\varphi}({\cal P})$,  take the corresponding
homogeneous binomial in variables $z_0, z_1, \ldots, z_n$
\[  B_{\varphi,z_0}({\cal P}) = z_{i_1}  \cdots
z_{i_k} - E_{\varphi}({\cal P})
z_{j_1}^{c_1}  \cdots z_{j_l}^{c_l}z_0^{(k - \sum_{s=1}^l c_s)}. \]
Then the elements $B_{\varphi,z_0}({\cal P})$ generate the ideal
$Q_{\varphi,z_0}(\Sigma)$, and K\"ahler limits of
\[ QH^*_{t \varphi}({\bf P}_{\Sigma}, {\bf C} \lbrack z_0, z_0^{-1} \rbrack),
\;\; t \rightarrow \infty  \]
 are
isomorphic to the Laurent polynomial extension
\[ H^*({\bf P}_{\Sigma}, {\bf C}) \lbrack z_0, z_0^{-1} \rbrack \]
of the odinary cohomology ring.
\end{theo}
\bigskip

Finally, if the first Chern class of ${\bf P}_{\Sigma}$ belongs to the
K\"ahler cone, i.e., $\alpha_{\Sigma} \in PL(\Sigma)$ is upper convex,
then it is possible to
define the quantum deformations of the cohomology ring of ${\bf P}_{\Sigma}$
over the polynomial ring ${\bf C}\lbrack z_0 \rbrack$.

\begin{opr}
{\rm  Assume that $\alpha_{\Sigma} \in PL(\Sigma)$ is upper convex.
We define the quantum cohomology ring
\[ QH^*_{\varphi}({\bf P}_{\Sigma}, {\bf C} \lbrack z_0 \rbrack) \]
over ${\bf C}\lbrack z_0 \rbrack$ as the quotion of the polynomial
ring  ${\bf C}\lbrack z_0, z_1, \ldots, z_n  \rbrack$ by the sum of
the ideal $P(\Sigma)\lbrack z_0 \rbrack $ and the ideal
\[ {\bf C}\lbrack z_0, z_1, \dots, z_n \rbrack \cap
Q_{\varphi,z_0}(\Sigma) \]
}
\end{opr}

\begin{theo}
 The ideal
\[ {\bf C}\lbrack z_0, z_1, \dots, z_n \rbrack \cap
Q_{\varphi,z_0}(\Sigma) \]
is generated by homogeneous binomials
\[   B_{\varphi,z_0}({\cal P}) = z_{i_1}  \cdots
z_{i_k} - E_{\varphi}({\cal P})
z_{j_1}^{c_1}  \cdots z_{j_l}^{c_l}z_0^{(k - \sum_{s=1}^l c_s)} \]
where ${\cal P}$ runs over all
primitive collections ${\cal P} \subset G(\Sigma)$.
$($Notice that convexity of $\alpha_{\Sigma}$
implies $k - \sum_{s=1}^l c_s \geq 0$.$)$

K\"ahler limits of the quantum cohomology ring
\[ QH^*_{t \varphi}({\bf P}_{\Sigma}, {\bf C} \lbrack z_0 \rbrack),
\;\; t \rightarrow \infty  \]
 are
isomorphic to the polynomial extension
\[ H^*({\bf P}_{\Sigma}, {\bf C}) \lbrack z_0  \rbrack \]
of the odinary cohomology ring.
\end{theo}

\section{Birational transformations}
\hspace*{\parindent}

It may look strange that we defined the  quantum cohomology rings
using infinitely many generators for the  ideals $Q_{\varphi}(\Sigma)$
and $Q_{\varphi, z_0}(\Sigma)$, while these ideals have only finite
number of generators indexed by primitive collections in $G(\Sigma)$.
The reason for that is the following important theorem:

\begin{theo}
Let $\Sigma_1$ and $\Sigma_2$ be two complete fans of regular cones such
that $G(\Sigma_1) = G(\Sigma_2)$, then the quantum cohomology rings
$QH^*_{\varphi}({\bf P}_{\Sigma_1}, {\bf C})$ and
$QH^*_{\varphi}({\bf P}_{\Sigma_2}, {\bf C})$
are isomorphic.
\label{flop}
\end{theo}

\Proof   Our definitions of  quantum cohomology rings does not
depend on the combinatiorial structure of the fan $\Sigma$, one needs to
know only all lattice vectors $v_1, \ldots, v_n \in G(\Sigma)$, but not the
combinatorial structure of the fan $\Sigma$. \hfill $\Box$
\bigskip

Since the equality $G(\Sigma_1) = G(\Sigma_2)$ means that two toric
varieties ${\bf P}_{\Sigma_1}$ and ${\bf P}_{\Sigma_2}$
are isomorphic in codimension $1$, we obtain

\begin{coro}
Let ${\bf P}_{\Sigma_1}$ and ${\bf P}_{\Sigma_2}$ be two smooth
compact toric manifolds which are
isomorphic in codimension $1$, then
the  rings
$QH^*_{\varphi}({\bf P}_1, {\bf C})$ and
$QH^*_{\varphi}({\bf P}_2, {\bf C})$
are isomorphic.
\end{coro}

\begin{exam}
{\rm Consider two $3$-dimensional fans $\Sigma_1$ and $\Sigma_2$ in
${\bf R}^3$ such that $G(\Sigma_1) = G(\Sigma_2) = \{ v_1, \ldots, v_6 \}$
where
\[ v_1 = (1,0,0),\;v_2 = (0,1,0),\;  v_3 = (0,0,1), \]
\[ v_4 = (-1,0,0),\;v_5 = (0,-1,0),\;  v_6 = (1,1,-1). \]
We define the combinatorial structure of  $\Sigma_1$ by
the primitive collections
\[ {\cal P}_1 = \{ v_1, v_4 \}, \;{\cal P}_2 = \{ v_2, v_5 \}, \;
{\cal P}_3 = \{ v_3, v_6 \}, \]
and the combinatorial structure of $\Sigma_2$ by
the primitive collections
\[ {\cal P}_1' = \{ v_1, v_4 \}, \;{\cal P}_2' = \{ v_2, v_5 \}, \;
{\cal P}_3' = \{ v_1, v_2 \}, \]
\[ {\cal P}_4' = \{ v_3, v_5, v_6 \}, \;{\cal P}_5' = \{ v_3, v_4, v_6 \}. \]

The flop between two toric manifolds is described by the diagrams:

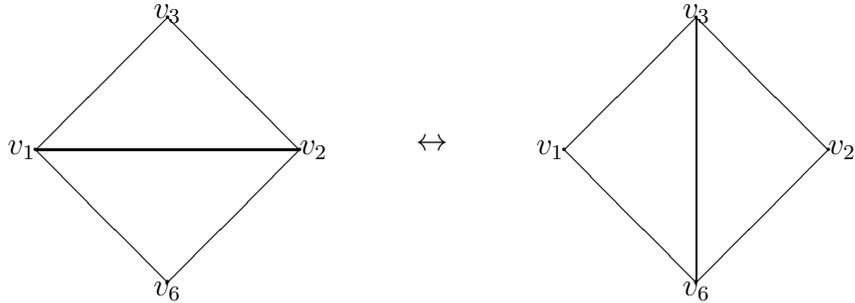
\begin{center}
\begin{picture}(360,200)
\put(100,150){\makebox(0,0){$\cdot$}}
\put(100,50){\makebox(0,0){$\cdot$}}
\put(50,100){\makebox(0,0){$\cdot$}}
\put(150,100){\makebox(0,0){$\cdot$}}

\put(300,150){\makebox(0,0){$\cdot$}}
\put(300,50){\makebox(0,0){$\cdot$}}
\put(250,100){\makebox(0,0){$\cdot$}}
\put(350,100){\makebox(0,0){$\cdot$}}

\put(100,150){\makebox(0,0)[b]{$v_3$}}
\put(100,50){\makebox(0,0)[t]{$v_6$}}
\put(50,100){\makebox(0,0)[r]{$v_1$}}
\put(150,100){\makebox(0,0)[l]{$v_2$}}

\put(300,150){\makebox(0,0)[b]{$v_3$}}
\put(300,50){\makebox(0,0)[t]{$v_6$}}
\put(250,100){\makebox(0,0)[r]{$v_1$}}
\put(350,100){\makebox(0,0)[l]{$v_2$}}

\put(200,100){\makebox(0,0)[b]{$\leftrightarrow$}}

\put(50,100){\line(1,1){50}}
\put(100,50){\line(1,1){50}}

\put(100,150){\line(1,-1){50}}
\put(50,100){\line(1,-1){50}}

\put(250,100){\line(1,1){50}}
\put(300,50){\line(1,1){50}}

\put(300,150){\line(1,-1){50}}
\put(250,100){\line(1,-1){50}}

\put(50,100){\line(1,0){100}}
\put(300,50){\line(0,1){100}}
\end{picture}
\end{center}

The ordinary cohomology rings $H^*({\bf P}_{\Sigma_1}, {\bf C})$ and
$H^*({\bf P}_{\Sigma_2}, {\bf C})$ are not isomorphic, because their
homogeneous ideals of
polynomial relations among $z_1, \ldots, z_6$   have different
numbers of minimal generators.
There exists the polynomial relation in the quantum cohomology
ring:
\[ \exp (\varphi(v_1) + \varphi(v_2)) z_1 z_2 =
\exp (\varphi(v_3) + \varphi(v_6)) z_3 z_6. \]
If $\varphi(v_1) + \varphi(v_2) <
\varphi(v_3) + \varphi(v_6)$, then we obtain the element
$z_3 z_6 \in SR(\Sigma_1)$ as  the limit for $t\varphi$, when
$t \rightarrow \infty$. On the other hand,  if
$\varphi(v_1) + \varphi(v_2) >   \varphi(v_3) + \varphi(v_6)$,
taking the same limit, we obtain $z_1 z_2 \in SR(\Sigma_2)$.
}
\end{exam}

Let us consider another simplest example of birational tranformation.

\begin{exam}
{\rm The quantum cohomology ring of the $2$-dimensional toric variety $F_1$
which is the blow-up of a point $p$ on ${\bf P}^2$
is isomorphic to the quotient
of the polynomial ring ${\bf C} \lbrack x_1, x_2
\rbrack$ by the ideal generated by two binomials
\[ x_1(x_1 + x_2) = \exp (-\phi_2);\; x_2^2 = \exp( - \phi_1) x_1, \]
where $x_1$ is the class of the $(-1)$-curve $C_1$ on $F_1$, $x_2$ is the
class of the fiber $C_2$ of the projection of $F_1$ on ${\bf P}^1$.
The numbers $\phi_1$ and $\phi_2$ are respectively
degrees of the restriction of the K\"ahler class
$\varphi$ on $C_1$ and $C_2$. }
\end{exam}

\begin{rem}
{\rm The definition of the quantum
cohomology ring for smooth toric manifolds
immediatelly extends to  the case of singular toric varieties. However,
the ordinary cohomology ring of singular toric varieties is not anymore
the K\"ahler limit of the quantum cohomology ring. In some cases, the
quantum cohomology ring of singular toric
varieties $V$ contains an information
about the ordinary cohomology ring of
special desingularizations $V'$ of $V$.
For instance, if we assume
that there exists a projective desingularization
$\psi\; : \; V' \rightarrow V$ such that $\psi^*{\cal K}_V = {\cal K}_{V'}$.
Then for every K\"ahler class $\varphi \in H^2 (V, {\bf C})$, one has
\[ {\rm dim}_{\bf C} H^*_{\varphi} (V, {\bf C}) =
{\rm dim}_{\bf C} H^*(V',{\bf C}). \]}
\end{rem}

\section{Geometric interpretation of quantum cohomology rings}
\hspace*{\parindent}

The spectra of the quantum cohomology ring
${\rm Spec}\,QH^*_{\varphi}({\bf P}_{\Sigma},
{\bf C} )$, and its two
polynomial versions
\[ {\rm Spec}\,QH^*_{\varphi}({\bf P}_{\Sigma},
{\bf C} \lbrack z_0, z_0^{-1} \rbrack),
\;  {\rm Spec}\,QH^*_{\varphi}({\bf P}_{\Sigma},
{\bf C} \lbrack z_0 \rbrack) \]
have simple geometric interpretations.

\begin{opr}
{\rm Denote by $\Pi(\Sigma)$ the $(n-d)$-dimensional affine subspace in
${\bf C}^n$ defined by the ideal $P(\Sigma)$.}
\end{opr}

\begin{opr}
{\rm Choose any isomorphism $N \cong {\bf Z}^d$, so that any element
$v \in N$ defines a Laurent monomial
$X^v$ in $d$ variables $X_1, \ldots, X_d$. Consider the embedding
of the $d$-dimensional torus $T(\Sigma) \cong ({\bf C}^*)^d$ in
$({\bf C}^*)^n$:
\[ (X_1, \ldots, X_d) \rightarrow (X^{v_1}, \ldots,  X^{v_n}). \]
Denote by $\Theta(\Sigma)$ the $(n-d)$-dimensional algebraic
torus $({\bf C}^*)^n /T(\Sigma)$.}
\end{opr}

\begin{opr}
{\rm Denote by $Exp$ the  analytical exponential mapping
\[ Exp\; : \; {\cal G} \rightarrow {\bf G} \]
where ${\bf G}$ is a complex analytic Lie group, and ${\cal G}$ is
its Lie algebra.

For example, one has the exponential mapping
\[ Exp \;:\; PL(\Sigma)_{\bf C} \rightarrow ({\bf C}^*)^n \]
\[ \varphi \mapsto
(e^{\varphi(v_1)}, \ldots, e^{\varphi(v_n)}) \]
which descends to the  exponential mapping
\[ Exp \; : \; H^2({\bf P}_{\Sigma}, {\bf C}) \rightarrow \Theta(\Sigma). \]
}
\end{opr}

\begin{prop}
The ${\bf T}(\Sigma)$-orbit ${T}_{\varphi}(\Sigma)$
of the point $Exp(\varphi) \in
({\bf C}^*)^n$ is closed, and its ideal is canonically isomorphic to
$Q_{\varphi}(\Sigma)$.
\end{prop}

\begin{coro}
The scheme ${\rm Spec}\,QH^*_{\varphi}({\bf P}_{\Sigma}, {\bf C} )$ is the
scheme-theoretic intersection of the  $d$-dimensional
subvariety $\overline{T}_{\varphi}(\Sigma) \subset
{\bf C}^n$ and the $(n-d)$-dimensional subspace $\Pi(\Sigma)$.
\end{coro}

\begin{opr}
{\rm Let $\tilde{N} = {\bf Z} \oplus N$. For any $v \in N$, define
$\tilde{v} \in \tilde{N}$ as $\tilde{v} = (1, v)$. Define
the embedding of the $(d+1)$-dimensional torus
$T^{\circ}(\Sigma) \cong ({\bf C}^*)^{d+1}$ in
$({\bf C}^*)^{n+1}$:
\[ (X_0, X_1, \ldots, X_d) \rightarrow
(X_0, X^{\tilde{v_1}}, \ldots,  X^{\tilde{v_n}}). \]}
\label{affine}
\end{opr}

The quotient $({\bf C}^*)^{n+1}/ T^{\circ}(\Sigma)$ is again isomorphic to
$\Theta(\Sigma)$.

\begin{prop}
The ideal of the ${T}^{\circ}(\Sigma)$-orbit
\[ T_{\varphi}^{\circ}(\Sigma) \subset
{\bf C}^* \times {\bf C}^n \]
of the  of the point
$(1, {Exp}\,(\varphi))  \in
({\bf C}^*)^{n+1}$ is canonically isomorphic to
$Q_{\varphi,z_0}(\Sigma)$.
\end{prop}

\begin{coro}
The scheme ${\rm Spect}\,QH^*_{\varphi}({\bf P}_{\Sigma},
{\bf C}\lbrack z_0, z_0^{-1}
\rbrack )$ is the
scheme-theo\-re\-tic intersection of the  $(d+1)$-dimensional
subvariety $T_{\varphi}^{\circ}(\Sigma) \subset
{\bf C}^* \times {\bf C}^{n}$ and the
$(n-d+1)$-dimensional subvariety
${\bf C}^* \times \Pi(\Sigma) \subset {\bf C}^* \times {\bf C}^n$.
\end{coro}

Similarly, one obtain the geometric interpretation of
$QH^*_{\varphi}({\bf P}_{\Sigma}, {\bf C}\lbrack z_0
\rbrack )$, when the first Chern class of ${\bf P}_{\Sigma}$ belongs to the
K\"ahler cone $K(\Sigma)$.

\begin{prop}
 The scheme ${\rm Spect}\,QH^*_{\varphi}({\bf P}_{\Sigma},
 {\bf C}\lbrack z_0
\rbrack )$ is the
scheme-theoretic intersection  in ${\bf C}^{n+1}$
of the  $(d+1)$-dimensional
${T}^{\circ}(\Sigma)$-orbit of the point $(1, Exp(\varphi))$ and
the
$(n-d+1)$-dimensional affine subspace ${\bf C} \times \Pi(\Sigma)
\subset {\bf C}^{n+1}$.
\end{prop}

The limits of quantum cohomology rings have also geometric interpretations.
One obtains, for instance, the spectrum of the ordinary cohomology
ring of ${\bf P}_{\Sigma}$
as the scheme-theoretic intersection of the affine subspace
$\Pi(\Sigma)$ with a "toric degeneration" of closures of
${\bf T}(\Sigma)$-orbits
$\overline{T}_{\varphi}(\Sigma) \subset
{\bf C}^n$. Such an interpretation allows to apply methods
of M. Kapranov, B. Strumfels, and A. Zelevinsky (see \cite{kap}, Theorem 5.3)
to establish connection between vertices of Chow polytope (secondary
polyhedron) and K\"ahler limits of quantum cohomology rings.

\section{Calabi-Yau hypersurfaces, Jacobian rings and the mirror symmetry}
\hspace*{\parindent}

Throughout in this section we fix a complete $d$-dimensional fan
of regular cones, and we assume that
${\bf P} = {\bf P}_{\Sigma}$ is a toric manifold
whose first Chern class belongs
to the closed K\"ahler cone $K(\Sigma)$, i.e., $\alpha: =
\alpha_{\Sigma}$ is a
convex ${\Sigma}$-piecewise linear function.

Let $\Delta = \Delta_{\alpha}$, the convex polyhedron in $M_{\bf R}$ (see
\ref{poly1}). For any sufficiently general section $S$ of
the anticanonical sheaf  ${\cal A}$ on ${\bf P}$ represented by
homogeneous polynomial $F(z)$,
the set $Z = \{ \pi(z) \in {\bf P} \; : \;  F(z) = 0  \}$ in ${\bf P}$
 is a
Calabi-Yau manifold ($c_1({\cal A}) = c_1({\bf P}))$.

Since the first Chern class of ${\bf P}$ in
the ordinary cohomology ring
$H^*({\bf P}, {\bf C})$ is the class of the
sum
$( z_1 + \cdots + z_n)$,
we obtain:

\begin{prop}
The image of $H^*({\bf P}, {\bf C})$ under the restriction
mapping to $H^*(Z, {\bf C})$ is isomorphic to the quotient
\[H^*({\bf P}, {\bf C})/ {\rm Ann}( z_1 + \cdots + z_n),  \]
where ${\rm Ann}( z_1 + \cdots + z_n)$ denotes the annulet of the
class of  $( z_1 + \cdots + z_n)$ in
$H^*({\bf P}, {\bf C})$.
\label{image}
\end{prop}

In general,  Proposition \ref{image} allows us to calculate only
a part of the ordinary cohomology
ring of a Calabi-Yau hypersurface $Z$ in toric variety
${\bf P}$. If
the first Chern class of ${\bf P}$ is in the interior of the
K\"ahler cone $K(\Sigma)$, then $Z$ is an ample divisor. For $d \geq 4$,
by Lefschetz theorem, the restriction mapping
$ H^2({\bf P}, {\bf C}) \rightarrow H^2(Z, {\bf C})$
is isomorphism. Thus, using Proposition \ref{image}, we can
calculate cup-products of any $(1,1)$-forms on $Z$.
\bigskip

\begin{opr}
{\rm Denote by
$\Delta^*$ the convex hull of the set $G(\Sigma)$ of all
generators, or equivalently,
\[ \Delta^* = \{ v \in N_{\bf R} \mid \alpha(v) \leq 1 \}. \]}
\end{opr}

\begin{rem}
{\rm The polyhedron $\Delta^*$ is dual to $\Delta$ reflexive polyhedron
(see \cite{bat.mir}). }
\end{rem}

\begin{theo}
There exists the canonical isomorphism between
the quantum cohomology ring
\[ QH^*_{\varphi}({\bf P}, {\bf C}  )   \]
and the Jacobian ring
\[ {\bf C} \lbrack X_1^{\pm 1}, \ldots, X_d^{\pm 1} \rbrack /
(X_1 \partial f/ \partial X_1, \ldots, X_d \partial f/ \partial X_d ) \]
of the Laurent polynomial
\[ f_{\varphi}(X) =
-1 + \sum_{i =1}^n \exp (\varphi(v_i))^{-1} X^{v_i}. \]
This isomorphism
is induced by the correspondence
\[ z_i \rightarrow X^{v_i}/\exp (\varphi(v_i)) \;\;
(1 \leq i \leq n). \]
In particular, it
maps the first Chern class $(z_1 + \ldots + z_n)$
of ${\bf P}$ to $f_{\varphi}(X) + 1 $.
\label{mirr}
\end{theo}

\Proof
Let
\[ {\cal H}\; : \; {\bf C} \lbrack z_1, \ldots, z_n \rbrack
\rightarrow {\bf C} \lbrack X_1^{\pm 1}, \ldots , X_d^{\pm 1} \rbrack \]
be the homomorphism defined by the correspondence
\[  z_i \rightarrow X^{v_i}/\exp (\varphi(v_i)).\]
By  \ref{def.fan}(iii),
${\cal H}$ is surjective. It is clear that
$Q_{\varphi}(\Sigma)$ is the kernel
of ${\cal H}$. On the other hand, if we a ${\bf Z}$-basis
$\{ u_1, \ldots, u_d \} \subset M$ which establishes isomorphisms
$M \cong {\bf Z}^d$ and $N \cong {\bf Z}^d$, we obtain:
\[{\cal H}(P(\Sigma)) = <
X_1 \partial f/ \partial X_1, \ldots, X_d \partial f/ \partial X_d >. \]
\hfill $\Box$
\bigskip

\begin{opr}
{\rm Let $S_{\Delta^*}$ be the affine coordinate ring of
the $T^{\circ}(\Sigma)$-orbit of the point
$(1,\ldots, 1) \in {\bf C}^{n+1}$
(see \ref{affine}). }
\end{opr}

\begin{opr}
{\rm For any Laurent polynomial
\[ f(X) = a_0 +  \sum_{i =1}^n  a_i X^{v_i}, \]
we define elements
\[ F_0, F_1, \ldots, F_d \in S_{\Delta^*} \]
as $F_i = \partial X_0 f(X) \partial X_0$, $(0 \leq i \leq d)$. }
\end{opr}

\begin{rem}
{\rm The ring $S_{\Delta^*}$ is a subring of
${\bf C}\lbrack X_0, X_1^{\pm 1}, \ldots, X_d^{\pm 1}
\rbrack$. There exists
the canonical grading of $S_{\Delta^*}$ by degree of $X_0$.

It is easy to see that the correspondence
\[ z_0 \rightarrow - X_0, \]
\[ z_i \rightarrow  X_0X^{v_i}/(\exp ( \varphi(v_i)))  \]
defines the isomorphism
\[ {\bf C} \lbrack z \rbrack / Q_{\varphi}(\Sigma) \cong
S_{\Delta^*}. \]
This isomorphism maps
$(- z_0 + z_1 + \cdots z_n )$ to $F_0$.}
\label{iso}
\end{rem}

\begin{theo}
{\rm (\cite{bat.var})}
Let
\[R_f = S_{\Delta^*} / < F_0, F_1, \ldots, F_d > .\]
Then the quotient
\[ R_f / {\rm Ann}\,(X_0) \]
is isomorphic to the $(d-1)$-weight subspace
$W_{d-1}H^{d-1}(Z_f, {\bf C})$ in
the cohomology space $H^{d-1}(Z_f, {\bf C})$
of the affine Calabi-Yau hypersurface
in $T(\Sigma)$ defined
by the Laurent polynomial $f(X)$.
\end{theo}

For any Laurent polynomial $f(X) = a_0 + \sum_{i=1}^n  a_i X^{v_i}$, we
can find an element $\varphi \in PL(\Sigma)_{\bf C}$ such that
\[ \frac{-a_i}{a_0} = \exp (- \varphi(v_i)). \]
A one-parameter family $t \varphi$ in $PL(\Sigma)$ induces
the one-parameter family of Laurent polynomials
\[ f_t (X) = - 1 +  \sum_{i =1}^n \exp (-  t \varphi(v_i))X^{v_i}. \]

Applying the isomorphism in \ref{iso}  and the statement in
Theorem \ref{basis}, we obtain the following:

\begin{theo}
Assume that $\varphi$ is in the interior of the K\"ahler cone $K(\Sigma)$.
Then the limit
\[ R_{f_t} / {\rm Ann } (X_0) \]
is isomorphic to
\[H^*({\bf P}, {\bf C})/ {\rm Ann}( z_1 + \cdots + z_n).  \]
\end{theo}

The last statement shows the relation, established in \cite{aspin4},
between the "toric" part of the topological cohomology rings of Calabi-Yau
$3$-folds in toric varieties and limits of the
multiplicative structure on $(d-1)$-weight part of the
Jacobian rings of their "mirrors".

\section{Topological sigma models on
toric manifolds}

So far we have not  explained why the ring $H^*_{\varphi} ({\bf P}_{\Sigma},
{\bf C})$ coincides
with the quantum cohomology ring corresponding to the topological
sigma model on $V$. In this section we want to establish the relations
between the ring $H^*_{\varphi} ({\bf P}_{\Sigma}, {\bf C})$ and
the quantum cohomology
rings  considered by physicists.
\bigskip

In order to apply the general
construction of the correlation functions in sigma models
(\cite{witten}, 3a ), we need the following information on the structure of
the space of holomorphic maps of
${\bf CP}^1$ to a $d$-dimensional toric manifold
${\bf P}_{\Sigma}$.

\begin{theo}
Let ${\cal I}$ be the moduli space of holomorphic maps $f\;:\; {\bf CP}^1
\rightarrow {\bf P}_{\Sigma}$. The space
${\cal I}$ consists of is infinitely many algebraic varieties
${\cal I}_{\lambda}$ indexed by elements
\[ \lambda = (\lambda_1, \ldots, \lambda_n) \in R(\Sigma), \]
where the numbers $\lambda_i$ are equal to the intersection
numbers ${\rm deg}_{{\bf CP}^1} f^*{\cal O}( Z_i)$ with divisors
$Z_i \subset {\bf P}_{\Sigma}$
such that $\pi^{-1}(Z_i)$ is defined by  the equation $z_i = 0$ in
$U(\Sigma)$.
Moreover, if all $\lambda_i \geq 0$,
then ${\cal I}_{\lambda}$ is irreducible
and the virtual dimension of ${\cal I}_{\lambda}$ equals
\[ d_{\lambda} = {\rm dim}_{\bf C} {\cal I}_{\lambda} =
d + \sum_{ i =1}^n \lambda_i. \]
\end{theo}

\Proof   The first statement follows immediatelly from the description
of the intersection product on ${\bf P}_{\Sigma}$ (\ref{deg}).

Assume now that all $\lambda_i$ are non-negative. This means the that the
preimage $f^{-1}(Z_i)$ consists of $\lambda_i$ points including their
multiplicities. Let ${\cal F}_{\Sigma}$ be the tangent bundle over
${\bf P}_{\Sigma}$. There exists the generalized Euler exact sequence
\[ 0 \rightarrow {\cal O}^{n-d}_{\bf P} \rightarrow {\cal O}_{\bf P}(Z_1)
\oplus \cdots \oplus {\cal O}_{\bf P}(Z_n) \rightarrow {\cal F}_{\Sigma}
\rightarrow 0. \]
Applying $f^*$, we obtain the short exact sequence of vector bundles on
${\bf CP}^1$.
\[ 0 \rightarrow {\cal O}^{n-d}_{{\bf CP}^1} \rightarrow
{\cal O}_{{\bf CP}^1}({\lambda}_1)
\oplus \cdots \oplus {\cal O}_{{\bf CP}^1}({\lambda}_n) \rightarrow
f^*{\cal F}_{\Sigma}
\rightarrow 0. \]
This implies that $h^1({\bf CP}^1,f^*{\cal F}_{\Sigma}) = 0$, and
$h^0({\bf CP}^1,f^*{\cal F}_{\Sigma}) = d + \lambda_1 + \cdots + \lambda_n$.

The irreducibility of ${\cal I}_{\lambda}$ for $\lambda \geq 0$ follows
from the explicit geometrical
construction of maps $f \in {\cal I}_{\lambda}$:

Choose $n$ polynomials $f_1(t), \ldots, f_n(t)$ such that
${\rm deg}\, f_i(t) = \lambda_i$ $( i =1, \ldots, n)$.
If all $\mid \lambda \mid =
\lambda_1 + \cdots + \lambda_n$ roots of $\{f_i\}$ are distinct, then
these polynomials define the mapping
\[ g\; :\; {\bf C} \rightarrow U(\Sigma) \subset {\bf C}^n .\]
The composition $\pi \circ g$ extends to the mapping $f$
of ${\bf CP}^1$ to ${\bf P}_{\Sigma}$ whose homology class is $\lambda$.
\hfill $\Box$
\bigskip

\begin{opr}
{\rm Let
\[ \Phi \; : \; {\cal I} \times {\bf CP}^1 \rightarrow {\bf P}_{\Sigma} \]
be the universal mapping. For every point $x \in {\bf CP}^1$ we denote by
$\Phi_x$ the restriction of $\Phi$ to ${\cal I}\times x$. The cohomology
classes $z_1 = \lbrack Z_1 \rbrack, \ldots, z_n = \lbrack Z_n \rbrack$ of
divisors $Z_1, \dots, Z_n$ on ${\bf P}_{\Sigma}$ in the ordinary cohomology
ring $H^*({\bf P}_{\Sigma})$ determine the cohomology
classes $W_{z_1}, \ldots, W_{z_n} \in H^*({\cal I})$ which are independent of
choice of $x \in {\bf CP}^1$. The element $W_{z_i}$ is the class of the
divisor on ${\cal I}$:
\[ \{ f \in {\cal I} \mid f(x) \in Z_i \}. \]
}
\end{opr}

The quantum cohomology ring of the sigma model with the target space
${\bf P}_{\Sigma}$ is defined by the intersection numbers
\[  (W_{\alpha_1} \cdot W_{\alpha_2} \cdot \cdots \cdot
W_{\alpha_k} )_{\cal I} \]
on the moduli space ${\cal I}$, where $\Phi_{\alpha} =
\Pi_x^*(\alpha)$

\begin{theo}
Let ${\bf P}_{\Sigma}$ be a $d$-dimensional toric manifold,
$\varphi \in H^2({\bf P}_{\Sigma}, {\bf C})$ a K\"ahler class.
Let $\lambda^0 = (\lambda_1^0, \ldots, \lambda_n^0)$ be a
non-negative element in $R(\Sigma)$,
$\Omega \in H^{2d}({\bf P}_{\Sigma}, {\bf C})$ the fundamental class
of the toric manifold ${\bf P}_{\Sigma}$. Then
the intersection number on the moduli space ${\cal I}$
\[ (W_{\Omega})\cdot
(W_{z_1})^{\lambda_1^0} \cdot(W_{z_2})^{\lambda_2} \cdot \cdots \cdot
(W_{z_n})^{\lambda_n^0} \]
vanishes
for all components ${\cal I}_{\lambda}$ except from $\lambda = \lambda_0$.
In the latter case, this number equals
\[ \exp (- {\rm deg}_{\varphi} \lambda). \]
\end{theo}

\Proof  Since the fundamental class $\Omega$ is involved
in the considered intersection number, this number  is zero for all
${\cal I}_{\lambda}$ such that the rational curves in the class $\lambda$
do not cover a dense Zariski open subset in ${\bf P}_{\Sigma}$. Thus,
we must consider only non-negative classes $\lambda$. Moreover,
the factors $(W_{z_i})^{\lambda_i^0}$ show that we must consider only those
$\lambda = (\lambda_1, \ldots, \lambda_n) \in R(\Sigma)$ such that
$\lambda_i \geq \lambda_i^0$, i.e., a mapping $f \in {\cal I}_{\lambda}$
is defined by polynomials $f_1, \ldots, f_n$ such that ${\rm deg}\, f_i
\geq \lambda_i$.

There is a general principle that non-zero contributions to
the intersection product
\[  (W_{\alpha_1} \cdot W_{\alpha_2} \cdot \cdots \cdot
W_{\alpha_k} )_{\cal I} \]
appear only from the components whose virtual ${\bf R}$-dimension
is equal to
\[ \sum_{i =1} {\rm deg} \, {\alpha_i}. \]
In our case, the last number is $d + \lambda_1^0 + \ldots + \lambda_n^0$.
Therefore,  a non-zero contribution
appears only if $\lambda = \lambda^0$.

It remains to notice that this contribution equals
$\exp(- {\rm deg}_{\varphi} \lambda_0)$.
The last statement follows from the observation that
the points $f^{-1}(Z_i) \subset {\bf CP}^1$ $( i = 1, \ldots, n)$
define the mapping $f\, : \,
{\bf CP}^1 \rightarrow {\bf P}_{\Sigma}$ uniquely up to the action of the
$d$-dimensional torus ${\bf T} = {\bf P}_{\Sigma} \setminus
(Z_1 \cup \cdots \cup Z_n)$, and the weight of the mapping $f$ in the
intersection  product is
\[ \int_{{\bf CP}^1} f^*(\varphi). \]\hfill $\Box$

\begin{coro}
Let ${\cal Z}_i$ be the quantum operator corresponding to the
class $\lbrack Z_i \rbrack \in H^2({\bf P}_{\Sigma}, {\bf C})$
$( i =1, \ldots, n)$ considered as an element
of the quantum cohomology ring. Then for every non-negative element
$\lambda \in R(\Sigma)$, one has the algebraic  relation
\[ {\cal Z}_1^{\lambda_1} \circ \cdots \circ {\cal Z}_n^{\lambda_n}  =
\exp (- {\rm deg}_{\varphi} \lambda)\, id. \]
\end{coro}

It turns out that the polynomial relations of above type are sufficient to
recover the quantum cohomology ring $H^*_{\varphi}({\bf P}_{\Sigma},
{\bf C})$:

\begin{theo}
Let $A_{\varphi}(\Sigma)$ be the quotient of the polynomial ring
${\bf C} \lbrack z \rbrack$ by the
sum of two ideals: $P(\Sigma)$ and the ideal generated by all
polynomials
\[ B_{\lambda} = z_1^{\lambda_1} \cdots z_n^{\lambda_n} -
\exp (-  {\rm deg}_{\varphi} \lambda) \]
where $\lambda$ runs over all non-negative elements of $R(\Sigma)$.
Then $A_{\varphi}(\Sigma)$ is isomorphic to
$H^*_{\varphi}({\bf P}_{\Sigma},
{\bf C})$.
\end{theo}

\Proof  Let $B_{\varphi}(\Sigma)$ be the ideal generated by all binomials
$B_{\lambda}$. By definition, $B_{\varphi}(\Sigma) \subset
Q_{\varphi}(\Sigma)$. So it is sufficient to prove that $Q_{\varphi}(\Sigma)
subset B_{\varphi}(\Sigma)$.

Let
\[ \sum_{i = 1}^n a_i v_i =  \sum_{j = 1}^n b_j v_j  \]
be a linear relation among $v_1, \ldots, v_n$ such that
$a_i, \, b_j \geq 0$. Since the set of
all nonnegative elements
$\lambda =(\lambda_1, \ldots, \lambda_n) \in R(\Sigma)$
$(\lambda_i \geq 0)$ generates a convex cone of maximal dimension in
$H^2({\bf P}_{\Sigma}, {\bf C})$, there exist two nonnegative
vectors $\lambda$, $\lambda' in R(\Sigma)$ such that
\[ \lambda - \lambda' = (\lambda_1 - \lambda_1', \ldots, \lambda_n -
\lambda_n') = (a_1 - b_1, \ldots, a_n - b_n ). \]
By definition, two binomials
$P_{\lambda}$ and $P_{\lambda'}$ are contained in $Q_{\varphi}(\Sigma)$.
Hence, the classes of $z_1, \ldots, z_n$ in ${\bf C} \lbrack z \rbrack/
B_{\varphi}(\Sigma)$ are invertible elements. Thus, the class of the binomial
\[ \exp ( \sum_{i = 1}^n a_i \varphi(v_i))
\prod_{i =1}^n z_j^{a_j} -
\exp ( \sum_{j = 1}^n b_j \varphi(v_j)) \prod_{j =1}^n z_j^{b_j} \]
is zero in ${\bf C} \lbrack z \rbrack/
B_{\varphi}(\Sigma)$. Thus,
$B_{\varphi}(\Sigma) = Q_{\varphi}(\Sigma)$. \hfill
$\Box$.

\end{document}